\documentclass[12pt,thmsa]{article}
\usepackage{amssymb}

\usepackage{sw20lart}

\input tcilatex
\begin{document}

\begin{center}
\textbf{Unfriendly or weakly unfriendly}

\textbf{partitions of graphs}

Francis Oger

29 janvier 2014\bigskip
\end{center}

\bigskip

\noindent \textbf{Abstract.} For each infinite cardinal $\kappa $\ and each
graph $G=(V,E)$, we say that a partition $\pi :V\rightarrow \left\{
0,1\right\} $\ is $\kappa $-unfriendly if, for each $x\in V$, $\left|
\left\{ y\in V\mid \left\{ x,y\right\} \in E\text{ and }\pi (y)\neq \pi
(x)\right\} \right| $ is $\geq \left| \left\{ y\in V\mid \left\{ x,y\right\}
\in E\text{ and }\pi (y)=\pi (x)\right\} \right| $ or $\geq \kappa $; $\pi $%
\ is unfriendly if the first property is true for each $x\in V$. Some
uncountable graphs of infinite minimum degree without unfriendly partition
have been constructed by S. Shelah and E.C. Milner, but it is not presently
known if\ countable graphs without unfriendly partition exist.

We show that, for each integer $n$, each graph of infinite minimum degree
has an $\omega _{n}$-unfriendly partition. We also prove that the following
properties are equivalent: (i) each graph has an $\omega $-unfriendly
partition; (ii) each countable graph has an unfriendly partition; (iii) each
countable graph without nonempty induced subgraph of infinite minimum degree
has an unfriendly partition (actually it is enough to consider a smaller
class of graphs).\bigskip

Here, a \emph{graph} is a pair $G=(V,E)$, where $V$ is the set of \emph{%
vertices} and $E$ is the set of \emph{edges} of $G$; the \emph{edges} are
non oriented pairs $\left\{ x,y\right\} $ with $x,y\in V$ and $x\neq y$. We
call\emph{\ induced graphs} the graphs $H=(W,F)$\ with $W\subset V$ and $%
F=\{\left\{ x,y\right\} \in E\mid x,y\in W\}$.

The \emph{neighbours} of a vertex $x$ are the vertices $y$ such that $%
\left\{ x,y\right\} \in E$. The \emph{degree} of $x$ is the cardinal of its
set of neighbours. The \emph{minimum degree} of $G$\ is the minimum of the
degrees of the vertices of $G$.

A \emph{partition} of $\Gamma $\ is a map $\pi :V\rightarrow \left\{
0,1\right\} $. We say that $\pi $ is \emph{unfriendly} if, for each $x\in V$%
, $\left| \left\{ y\in V\mid \left\{ x,y\right\} \in E\text{ and }\pi
(y)\neq \pi (x)\right\} \right| $ is

\noindent $\geq \left| \left\{ y\in V\mid \left\{ x,y\right\} \in E\text{
and }\pi (y)=\pi (x)\right\} \right| $ (the two sets can be infinite).

\noindent For each cardinal $\kappa $, we say that $\pi $\ is $\kappa -\emph{%
unfriendly}$ if, for each $x\in V$,

\noindent $\left| \left\{ y\in V\mid \left\{ x,y\right\} \in E\text{ and }%
\pi (y)\neq \pi (x)\right\} \right| $ is $\geq \kappa $ or

\noindent $\geq \left| \left\{ y\in V\mid \left\{ x,y\right\} \in E\text{
and }\pi (y)=\pi (x)\right\} \right| $.\bigskip

\textbf{1. Graphs of infinite minimum degree.}\bigskip

Here we consider the following question: For which cardinals $\kappa $\ is
it true that any graph of infinite minimum degree has a $\kappa $-unfriendly
partition?\bigskip

\noindent \textbf{Theorem 1.1.} Let $\kappa $\ be an infinite cardinal and
let $G=(V,E)$\ be a graph of infinite minimum degree. Suppose that each
induced graph $H=(W,F)$ with $\left| W\right| \leq \kappa $ has an
unfriendly partition if each element of $W$ has $\kappa $\ neighbours in $W$%
, or more neighbours in $W$ than in $V-W$. Then $G$ has a $\kappa $%
-unfriendly partition.\bigskip

\noindent \textbf{Proof.} Let $\mathcal{E}$ consist of the pairs $(W,\pi )$,
with $W\subset V$\textbf{\ }and $\pi :W\rightarrow \left\{ 0,1\right\} $,
such that:

\noindent each element of $W$ has at least $\kappa $\ neighbours in $W$ or
at least as many neighbours in $W$ as in $V-W$;

\noindent $\pi $ is a $\kappa $-unfriendly partition of the induced graph
defined on $W$.

Then $\mathcal{E}$ contains $(\varnothing ,\rho )$ where $\rho :\varnothing
\rightarrow \left\{ 0,1\right\} $\ is the trivial map. For each cardinal $%
\mu $, the union of any increasing sequence $(W_{\alpha },\pi _{\alpha
})_{\alpha <\mu }$ of elements of $\mathcal{E}$ belongs to $\mathcal{E}$.

Consequently, it suffices to prove that, for each $(W,\pi )\in \mathcal{E}$
with $W\subsetneq V$, there exists\ $(W^{\prime },\pi ^{\prime })\in 
\mathcal{E}$\ with $W\subsetneq W^{\prime }$.

First suppose that there exists $x\in V-W$\ with at least $\kappa $
neighbours in $W$, or at least as many neighbours in $W$\ as in $V-W$.\ Then
we write $W^{\prime }=W\cup \left\{ x\right\} $\textbf{\ }and $\pi ^{\prime
}(x)=1$\ except if

\noindent $\left| \left\{ y\in W\mid \left\{ x,y\right\} \in E\text{\textbf{%
\ }and }\pi (y)=1\right\} \right| >\left| \left\{ y\in W\mid \left\{
x,y\right\} \in E\text{\textbf{\ }and }\pi (y)=0\right\} \right| $.

Now we can suppose that each $y\in V-W$ has less than $\kappa $ neighbours
in $W$\textbf{\ }and more neighbours in $V-W$\ than in $W$. Then we consider
a vertex $x\in V-W$\textbf{\ }and we define by induction an increasing
sequence $(A_{n})_{n\in 
%TCIMACRO{\U{2115} }%
%BeginExpansion
\mathbb{N}
%EndExpansion
}$ of subsets of $V-W$. We write $A_{0}=\left\{ x\right\} $. For each $n\in 
%TCIMACRO{\U{2115} }%
%BeginExpansion
\mathbb{N}
%EndExpansion
$, $A_{n+1}$ is obtained from $A_{n}$ by adding, for each $y\in A_{n}$:

\noindent the neighbours of $y$ in $V-W$ if there exist at most $\kappa $ of
them;

\noindent $\kappa $ neighbours of $y$ in $V-W$ if there exist more than $%
\kappa $ of them.

We write $A=\cup _{n\in 
%TCIMACRO{\U{2115} }%
%BeginExpansion
\mathbb{N}
%EndExpansion
}A_{n}$. We have $\left| A\right| \leq \kappa $. For each $n\in 
%TCIMACRO{\U{2115} }%
%BeginExpansion
\mathbb{N}
%EndExpansion
$\textbf{\ }and each $y\in A_{n}$, if $y$ has at least $\kappa $\ neighbours
in $V$,\textbf{\ }and therefore in $V-W$, then it has\ $\kappa $\ neighbours
in $A_{n+1}$\textbf{\ }and therefore in $A$. Otherwise, $y$ has more
neighbours in $A_{n+1}$\ than in $V-A_{n+1}$ since it has more neighbours in 
$V-W$\ than in $W$; consequently, $y$ has more neighbours in $A$\ than in $%
V-A$.

Now, according to the hypotheses of the Theorem, the induced graph defined
on $A$ admits an unfriendly partition $\tau $. It follows that $\pi ^{\prime
}=\tau \cup \pi $\ is a $\kappa $-unfriendly partition of the induced graph
defined on $A\cup W$.~~$\blacksquare $\bigskip 

By [2, Theorem 1], there exists a graph $G=(V,E)$ of minimum degree $\omega $
with $\left| V\right| =(2^{\omega })^{(+\omega )}$\ which has no unfriendly
partition. Moreover, by [2, Theorem 3], it is consistant that there exists a
graph $G=(V,E)$ of minimum degree $\omega $ with $\left| V\right| =\omega
_{\omega }$\ which has no unfriendly partition. Consequently, in the case of
graphs of infinite minimum degree, the following result is the best possible
one:\bigskip

\noindent \textbf{Corollary 1.2.} For each $n\in 
%TCIMACRO{\U{2115} }%
%BeginExpansion
\mathbb{N}
%EndExpansion
$, each graph of infinite minimum degree admits an $\omega _{n}$-unfriendly
partition.\bigskip

\noindent \textbf{Proof.} Il follows from [1, Theorem 2] that, for each $%
n\in 
%TCIMACRO{\U{2115} }%
%BeginExpansion
\mathbb{N}
%EndExpansion
$, each graph of cardinal $\omega _{n}$ with infinite minimum degree admits
an unfriendly partition. The Corollary is a consequence of Theorem 1.1 and
this fact.~~$\blacksquare $\bigskip

\textbf{2. Graphs of arbitrary minimum degree.}\bigskip

Now we consider the following question: Does each graph $G=(V,E)$ admit an $%
\omega $-unfriendly partition? The answer is not presently known, even for
countable graphs. We show that a positive answer for a particular class of
countable graphs implies a positive answer for all graphs.\bigskip 

\noindent \textbf{Definitions.} A \emph{graph with finite conditions} (resp. 
\emph{with conditions}) is a pair $\Gamma =(G,\mathcal{C})$\ where $G=(V,E)$
is a graph\textbf{\ }and $\mathcal{C}=(\kappa _{x},\lambda _{x})_{x\in V}$\
with $\kappa _{x},\lambda _{x}$ finite (resp. finite or infinite) cardinals
for each $x\in V$.

We say that a map $\pi :V\rightarrow \left\{ 0,1\right\} $ is an \emph{%
unfriendly} \emph{partition }of $\Gamma $ if, for each $x\in V$, $\pi (x)=0$
implies

\noindent $\kappa _{x}+\left| \left\{ y\in V\mid \left\{ x,y\right\} \in E%
\text{\textbf{\ }and }\pi (y)=0\right\} \right| \leq \lambda _{x}+\left|
\left\{ y\in V\mid \left\{ x,y\right\} \in E\text{\textbf{\ }and }\pi
(y)=1\right\} \right| $

\noindent and $\pi (x)=1$ implies

\noindent $\kappa _{x}+\left| \left\{ y\in V\mid \left\{ x,y\right\} \in E%
\text{\textbf{\ }and }\pi (y)=0\right\} \right| \geq \lambda _{x}+\left|
\left\{ y\in V\mid \left\{ x,y\right\} \in E\text{\textbf{\ }and }\pi
(y)=1\right\} \right| $.

\noindent For each cardinal $\mu $, we say that $\pi $ is a\emph{\ }$\mu $%
\emph{-unfriendly} \emph{partition} of $\Gamma $ if the two implications
above are true, except possibly when the two cardinals considered are $\geq
\mu $.\bigskip

\noindent \textbf{Proposition 2.1.} Let $\Gamma =(G,\mathcal{C})$\ be a
graph with conditions. Then there exists a set $\mathcal{C}^{\prime }$ of
finite conditions, defined on $G$\textbf{\ }and trivial on the elements of
infinite degree, such that each $\omega $-unfriendly partition of $\Gamma
^{\prime }=(G,\mathcal{C}^{\prime })$ is an $\omega $-unfriendly partition
of $\Gamma $.\bigskip

\noindent \textbf{Proof.} We write $G=(V,E)$\textbf{\ }and $\mathcal{C}%
=(\kappa _{x},\lambda _{x})_{x\in V}$. The set $\mathcal{C}^{\prime
}=(\kappa _{x}^{\prime },\lambda _{x}^{\prime })_{x\in V}$ is defined as
follows:

\noindent $\kappa _{x}^{\prime }=\lambda _{x}^{\prime }=0$ for $x$ of
infinite degree;

\noindent $\kappa _{x}^{\prime }=\kappa _{x}$\textbf{\ }and $\lambda
_{x}^{\prime }=\lambda _{x}$ for $x$ of finite degree\textbf{\ }and $\kappa
_{x},\lambda _{x}$\ finite;

\noindent $\kappa _{x}^{\prime }=n+1$\textbf{\ }and $\lambda _{x}^{\prime
}=0 $ for $x$ of finite degree $n$\textbf{\ }and $\kappa _{x}$\ infinite;

\noindent $\kappa _{x}^{\prime }=0$\textbf{\ }and $\lambda _{x}^{\prime
}=n+1 $ for $x$ of finite degree $n$, $\kappa _{x}$\ finite\textbf{\ }and $%
\lambda _{x}$\ infinite.~~$\blacksquare $\bigskip

\noindent \textbf{Notation.} For each graph $G=(V,E)$, we denote by $%
\mathcal{I}(G)$ the largest $W\subset V$ such that each element of $W$ has
infinitely many neighbours in $W$.\bigskip

\noindent \textbf{Remark.} We have $\mathcal{I}(G)=V$ if\textbf{\ }and only
if the minimum degree of $G$ is infinite. We can have $\mathcal{I}%
(G)=\emptyset $.\bigskip 

\noindent \textbf{Theorem 2.2.} Let $\Gamma =(G,\mathcal{C})$\ be a graph
with conditions. Suppose that each countable graph with finite conditions $%
\Delta =(H,\mathcal{D})$, with $H=(W,F)$ induced by $G$\textbf{\ }and $W\cap 
\mathcal{I}(G)=\emptyset $, has an unfriendly partition. Then $\Gamma $ has
an $\omega $-unfriendly partition.\bigskip

\noindent \textbf{Proof.} We write $G=(V,E)$\textbf{\ }and $\mathcal{C}%
=(\kappa _{x},\lambda _{x})_{x\in V}$. We define by induction on the ordinal 
$\alpha $ two sequences $(V_{\alpha })_{\alpha \leq \delta }$\textbf{\ }and $%
(W_{\alpha })_{\alpha <\delta }$\ of subsets of $V$.

We write $V_{0}=\emptyset $\textbf{\ }and $W_{0}=\mathcal{I}(G)$. For each $%
\alpha \geq 1$, supposing $V_{\beta }$\textbf{\ }and $W_{\beta }$\ already
defined for each $\beta <\alpha $, we write $V_{\alpha }=\cup _{\beta
<\alpha }W_{\beta }$; if $V_{\alpha }=V$, then we write $\delta =\alpha $%
\textbf{\ }and we do not define\ $W_{\alpha }$; otherwise, we take for\ $%
W_{\alpha }$ any countable subset $X$ of $V-V_{\alpha }$\ such that, for
each $x\in X$, all the neighbours of $x$ in $V-V_{\alpha }$, or infinitely
many of them, belong to $X$.

For each $\alpha <\delta $, we consider the induced graph $H_{\alpha }$
defined on $W_{\alpha }$. We define by induction on $\alpha <\delta $ a
sequence of conditions $\mathcal{D}_{\alpha }=(\kappa _{\alpha ,x},\lambda
_{\alpha ,x})_{x\in W_{\alpha }}$\textbf{\ }and a map $\pi _{\alpha
}:W_{\alpha }\rightarrow \left\{ 0,1\right\} $.

We write $\mathcal{D}_{0}=(\kappa _{x},\lambda _{x})_{x\in W_{0}}$. As $%
H_{0} $ is a graph of infinite minimum degree, it admits an $\omega $%
-unfriendly partition $\pi _{0}$ by Theorem 1.1. Then $\pi _{0}$ is also an $%
\omega $-unfriendly partition of $(H_{0},\mathcal{D}_{0})$.

For $1\leq \alpha <\delta $, supposing $\pi _{\beta }$\textbf{\ }and $%
\mathcal{D}_{\beta }$\ already defined for each $\beta <\alpha $, we write $%
\kappa _{\alpha ,x}=\kappa _{x}+\left| \cup _{\beta <\alpha }\left\{ y\in
W_{\beta }\mid \left\{ x,y\right\} \in E\text{\textbf{\ }and }\pi _{\beta
}(y)=0\right\} \right| $\ and $\lambda _{\alpha ,x}=\lambda _{x}+\left| \cup
_{\beta <\alpha }\left\{ y\in W_{\beta }\mid \left\{ x,y\right\} \in E\text{%
\textbf{\ }and }\pi _{\beta }(y)=1\right\} \right| $ for each $x\in
W_{\alpha }$. It follows from Proposition 2.1 and the hypotheses\textbf{\ }%
of the Theorem that $(H_{\alpha },\mathcal{D}_{\alpha })$\ admits an $\omega 
$-unfriendly partition $\pi _{\alpha }$.

The map $\pi =\cup _{\alpha <\delta }\pi _{\alpha }$\ is an $\omega $%
-unfriendly partition of $\Gamma $.~~$\blacksquare $\bigskip

Now we show that, in order to prove the existence of unfriendly partitions
for countable graphs with conditions, it suffices to consider countable
graphs without conditions. By Proposition 2.1, it is enough to consider the
pairs $(G,\mathcal{C})$ where $G=(V,E)$ is a countable graph and $\mathcal{C}%
=(m_{x},n_{x})_{x\in V}$ is a set of finite conditions with $m_{x}=n_{x}=0$\
for each $x$\ of infinite degree.

For each such pair, we define a graph $H=(W,F)$ as follows: We consider the
union of three disjoint copies $((W_{1},F_{1}),\mathcal{D}_{1})$, $%
((W_{2},F_{2}),\mathcal{D}_{2})$, $((W_{3},F_{3}),\mathcal{D}_{3})$ of $(G,%
\mathcal{C})$. We introduce some new vertices:

\noindent $u_{1},u_{2},u_{3}$, $(v_{i,n})_{i=1,2,3;n\in 
%TCIMACRO{\U{2115} }%
%BeginExpansion
\mathbb{N}
%EndExpansion
}$,\ $(w_{i,n})_{i=1,2,3;n\in 
%TCIMACRO{\U{2115} }%
%BeginExpansion
\mathbb{N}
%EndExpansion
}$,\ $(x_{i,n})_{i=1,2,3;n\in 
%TCIMACRO{\U{2115} }%
%BeginExpansion
\mathbb{N}
%EndExpansion
}$,

\noindent and some new edges:

\noindent $\left\{ u_{i},v_{j,n}\right\} $\ for $i\neq j$ in $\{1,2,3\}$\
and $n\in 
%TCIMACRO{\U{2115} }%
%BeginExpansion
\mathbb{N}
%EndExpansion
$;

\noindent $\left\{ v_{i,3n},w_{i,n}\right\} $, $\left\{
v_{i,3n+1},w_{i,n}\right\} $, $\left\{ v_{i,3n+2},w_{i,n}\right\} $\ for $%
i\in \{1,2,3\}$\ and $n\in 
%TCIMACRO{\U{2115} }%
%BeginExpansion
\mathbb{N}
%EndExpansion
$;

\noindent $\left\{ w_{i,2n},x_{i,n}\right\} $, $\left\{
w_{i,2n+1},x_{i,n}\right\} $\ for $i\in \{1,2,3\}$\ and $n\in 
%TCIMACRO{\U{2115} }%
%BeginExpansion
\mathbb{N}
%EndExpansion
$.

For each $i\in \left\{ 1,2,3\right\} $, as a substitute to $\mathcal{D}_{i}$%
, we put $m_{y}$\ edges of type\ $\left\{ w_{i,n},y\right\} $ and $n_{y}$\
edges of type\ $\left\{ x_{i,n},y\right\} $ for each $y\in W_{i}$, all of
them defined in such a way that each $w_{i,n}$, and each $x_{i,n}$, is an
endpoint of at most one edge of that type.

Now we show that each unfriendly partition $\pi $\ of $H$ induces an
unfriendly partition of $(G,\mathcal{C})$. We consider $i\neq j$ in $%
\{1,2,3\}$\ such that $\pi (u_{i})=\pi (u_{j})$. We can suppose for instance
that $i=2$ and $j=3$. Replacing if necessary $\pi $ by $1-\pi $, we can also
suppose $\pi (u_{2})=\pi (u_{3})=0$.

Then we necessarily have $\pi (v_{1,n})=1$, $\pi (w_{1,n})=0$\ and $\pi
(x_{1,n})=1$\ for each $n\in 
%TCIMACRO{\U{2115} }%
%BeginExpansion
\mathbb{N}
%EndExpansion
$. It follows that $\pi $ induces an unfriendly partition of $((W_{1},F_{1}),%
\mathcal{D}_{1})$, and therefore an unfriendly partition of $(G,\mathcal{C})$%
.

We observe that, in $H$, the only vertices of infinite degree are:

\noindent $u_{1},u_{2},u_{3}$, which have no neighbour of infinite degree;

\noindent for $i\in \left\{ 1,2,3\right\} $, the vertices of infinite degree
of $(W_{i},F_{i})$, whose only neighbours in $H$\ are their neighbours in\ $%
(W_{i},F_{i})$.

\noindent In particular, $\mathcal{I}(G)=\emptyset $ implies $\mathcal{I}%
(H)=\emptyset $.

In view of Theorem 2.2, we have:\bigskip

\noindent \textbf{Corollaire 2.3.} The following properties are equivalent:

\noindent 1) each graph with conditions has an $\omega $-unfriendly
partition;

\noindent 2) each graph has an $\omega $-unfriendly partition;

\noindent 3) each countable graph has an unfriendly partition;

\noindent 4) each countable graph without nonempty induced subgraph of
infinite minimum degree has an unfriendly partition.\bigskip

The proposition below implies that we can consider an even smaller class of
graphs in order to prove the existence of unfriendly partitions:\bigskip

\noindent \textbf{Proposition 2.4.} For each graph $G=(V,E)$, there exists a
unique partition $V=V^{1}\cup V^{2}\cup \mathcal{I}(G)$\ such that each
element of $V^{1}$ has finitely many neighbours in $V^{1}$ and each element
of $V^{2}$ has infinitely many neighbours in $V^{1}$.\bigskip

\noindent \textbf{Proof.} We define by induction on the ordinal $\alpha $
the subsets $V_{\alpha }^{1},V_{\alpha }^{2}$ with

\noindent $V_{\alpha }^{1}=\{x\in V-\cup _{\beta <\alpha }V_{\beta
}^{1}-\cup _{\beta <\alpha }V_{\beta }^{2}\mid x$ has finitely many
neighbours in $V-\cup _{\beta <\alpha }V_{\beta }^{2}\}$ and

\noindent $V_{\alpha }^{2}=\{x\in V-\cup _{\beta \leq \alpha }V_{\beta
}^{1}-\cup _{\beta <\alpha }V_{\beta }^{2}\mid x$ has finitely many
neighbours in $V-\cup _{\beta \leq \alpha }V_{\beta }^{1}-\cup _{\beta
<\alpha }V_{\beta }^{2}\}$.

\noindent We write $V^{1}=\cup _{\gamma <\delta }V_{\gamma }^{1}$, $%
V^{2}=\cup _{\gamma <\delta }V_{\gamma }^{2}$ and $V^{3}=V-(V^{1}\cup V^{2})$
where $\delta $ is the smallest integer such that $V_{\delta }^{1}=V_{\delta
}^{2}=\emptyset $.

It follows from the definition of the subsets $V_{\alpha }^{1}$\ that each
element of $V^{1}$\ has finitely many neighbours in $V^{1}\cup V^{3}$.

For each $\alpha <\delta $, each $x\in V_{\alpha }^{2}$\ has infinitely many
neighbours in $V-\cup _{\beta <\alpha }V_{\beta }^{2}$\ since it does not
belong to $V_{\alpha }^{1}$, and therefore infinitely many neighbours in $%
\cup _{\beta \leq \alpha }V_{\beta }^{1}$\ since it only has finitely many
neighbours in $V-\cup _{\beta \leq \alpha }V_{\beta }^{1}-\cup _{\beta
<\alpha }V_{\beta }^{2}$. It follows that each element of $V^{2}$ has
infinitely many neighbours in $V^{1}$.

The definition of $\delta $\ implies that each element of $V^{3}$ has
infinitely many neighbours in $V^{3}$. Consequently, we have $V^{3}\subset 
\mathcal{I}(G)$.

It remains to be proved that $\mathcal{I}(G)\subset V^{3}$. We show by
induction on the ordinal $\alpha $ that $V_{\alpha }^{1}\cap \mathcal{I}%
(G)=\emptyset $ and\ $V_{\alpha }^{2}\cap \mathcal{I}(G)=\emptyset $. If
these two properties are true for each $\beta <\alpha $, then we have $%
V_{\alpha }^{1}\cap \mathcal{I}(G)=\emptyset $ since each $x\in V_{\alpha
}^{1}$ has finitely many neighbours in $V-\cup _{\beta <\alpha }V_{\beta
}^{2}$, which contains $\mathcal{I}(G)$ by the induction hypothesis.\
Similarly, we have $V_{\alpha }^{2}\cap \mathcal{I}(G)=\emptyset $ because
each $x\in V_{\alpha }^{2}$ has finitely many neighbours in $V-\cup _{\beta
\leq \alpha }V_{\beta }^{1}-\cup _{\beta <\alpha }V_{\beta }^{2}$, which
contains $\mathcal{I}(G)$ by the induction hypothesis since $V_{\alpha
}^{1}\cap \mathcal{I}(G)=\emptyset $.

Now let us consider another partition $V=W^{1}\cup W^{2}\cup \mathcal{I}(G)$%
\ such that each element of $W^{1}$ has finitely many neighbours in $W^{1}$
and each element of $W^{2}$ has infinitely many neighbours in $W^{1}$. Then
we prove by induction on the ordinal $\alpha $ that $V_{\alpha }^{1}\subset
W^{1}$\ and $V_{\alpha }^{2}\subset W^{2}$. If these two properties are true
for each $\beta <\alpha $, then we have $V_{\alpha }^{1}\subset W^{1}$ since
each $x\in V_{\alpha }^{1}$ has finitely many neighbours in $V-\cup _{\beta
<\alpha }V_{\beta }^{2}$, which contains $W^{1}$ by the induction
hypothesis.\ Similarly, we have $V_{\alpha }^{2}\subset W^{2}$ because each $%
x\in V_{\alpha }^{2}$ has infinitely many neighbours in $\cup _{\beta \leq
\alpha }V_{\beta }^{1}$, which is contained in $W^{1}$ by the induction
hypothesis since $V_{\alpha }^{1}\subset W^{1}$.~~$\blacksquare $\bigskip

By Theorem 2.2, we can suppose that the graph $G=(V,E)$ of Proposition 2.4
satisfies $\mathcal{I}(G)=\emptyset $. Then we consider the set $V^{F}$
consisting of its vertices of finite degree, and the graph $G^{\ast
}=(V,E^{\ast })$, where $E^{\ast }$\ is obtained from $E$ by deleting the
vertices $\left\{ x,y\right\} $\ for $x,y\in V^{1}-V^{F}$ and for $x,y\in
V^{2}$. Any unfriendly partition of $G^{\ast }$ is also an unfriendly
partition of $G$.

Consequently, it suffices to prove the existence of unfriendly partitions
for the countable graphs $G=(V,E)$ such that $V=V^{1}\cup V^{2}$ and such
that there exists no edge $\left\{ x,y\right\} $\ for $x,y\in V^{1}-V^{F}$
and for $x,y\in V^{2}$.

In view of this simplification, we propose the question below. A positive
answer would be a first step to prove that each countable graph has an
unfriendly partition. We note that [3] does not give an answer to this
question, even though it proves the existence of unfriendly partitions for
large classes of countable graphs.\bigskip

\noindent \textbf{Question.} Does any countable graph $G=(V,E)$ admit an
unfriendly partition if each element of $V-V^{F}$ has infinitely many
neighbours in $V^{F}$ and no neighbour in $V-V^{F}$?

\bigskip

\begin{center}
\textbf{References}\bigskip
\end{center}

\noindent \lbrack 1] R. Aharoni, E.C. Milner and K. Prikry, Unfriendly
partitions of a graph, J. Combinatorial Theory (series B) 50 (1990), 1-10.

\noindent \lbrack 2] S. Shelah and E.C. Milner, Graphs with no unfriendly
partition, in A tribute to Paul Erd\"{o}s (pp. 373-384), Cambridge
University Press, Cambridge, 1990.

\noindent \lbrack 3] H. Bruhn, R. Diestel, A. Georgakopoulos and P. Spr\"{u}%
ssel, Every rayless graph has an unfriendly partition, Combinatoria 30
(2010), 521-532.

\vfill

\end{document}